\documentclass[12pt,a4paper]{amsart}
\usepackage{amssymb,amsmath}

\headsep=1cm \numberwithin{equation}{section}
\newtheorem{theorem}{Theorem}[section]
\newtheorem{lemma}[theorem]{Lemma}

\theoremstyle{definition}
\newtheorem{definition}[theorem]{Definition}
\theoremstyle{remark}

\newtheorem{question}[theorem]{Question}

\newcommand{\Ass}{\operatorname{Ass}}
\newcommand{\im}{\operatorname{im}}
\newcommand{\coker}{\operatorname{coker}}

\newcommand{\pd}{\operatorname{pd}}
\newcommand{\V}{\operatorname{V}}

\newcommand{\Ext}{\operatorname{Ext}}
\newcommand{\Supp}{\operatorname{Supp}}
\newcommand{\Tor}{\operatorname{Tor}}
\newcommand{\Hom}{\operatorname{Hom}}

\newcommand{\BN}{\Bbb N}

\newcommand{\lo}{\longrightarrow}
\newcommand{\fm}{\frak{m}}
\newcommand{\fp}{\frak{p}}

\newcommand{\fa}{\frak{a}}

\def\mapdown#1{\Big\downarrow\rlap
{$\vcenter{\hbox{$\scriptstyle#1$}}$}}

\begin{document}
\author[Divaani-Aazar and Sazeedeh]{Kamran Divaani-Aazar and Reza Sazeedeh}
\title[Cofiniteness of generalized local cohomology modules]
{Cofiniteness of generalized local cohomology modules}

\address{K. Divaani-Aazar, Department of Mathematics, Az-Zahra University,
Vanak, Post Code 19834, Tehran, Iran and Institute for Studies in
Theoretical Physics and Mathematics, P.O. Box 19395-5746, Tehran,
Iran.} \email{kdivaani@ipm.ir}

\address{R. Sazeedeh, Department of Mathematics, Uromeiyeh University, Uromeiyeh,
Iran and Institute of Mathematics, University for Teacher
Education, 599 Taleghani Avenue, Tehran 15614, Iran.}



\begin{abstract}
Let $\fa$ denote an ideal of a commutative Noetherian ring $R$ and
$M$ and $N$ two finitely generated $R$-modules with $\pd M<
\infty$. It is shown that if $\fa$ is principal or $R$ is complete
local and $\fa$ a prime ideal with $\dim R/\fa=1$, then the
generalized local cohomology module $H^i_{\fa}(M,N)$ is
$\fa$-cofinite for all $i \geq 0$. This provides an affirmative
answer for the above ideal $\fa$  to a question proposed in [{\bf
13}].
\end{abstract}

\maketitle

\section{Introduction}

A generalization of local cohomology functors has been given by J.
Herzog in [{\bf 6}]. Let $\fa$ denote an ideal of a commutative
Noetherian ring $R$. For each $i\geq 0$, the functor
$H^i_{\fa}(.,.)$ defined by
$H^i_{\fa}(M,N)=\underset{n}{\varinjlim}\Ext_R^i(M/\fa^nM,N)$, for
all $R$-modules $M$ and $N$. Clearly, this notion is a
generalization of the usual local cohomology functor. The study of
this concept was continued in the articles [{\bf 10}], [{\bf 2}]
and [{\bf 12}]. Recently, there is some new interest on studying
generalized local cohomology (see e.g. [{\bf 1}], [{\bf 13}] and
[{\bf 14}]).

In 1969, A. Grothendieck conjectured that if $\fa$ is an ideal of
$R$ and $N$ is a finitely generated $R$-module, then
$\Hom_R(R/\fa,H^i_{\fa}(N))$ is finitely generated for all $i\geq
0$. R. Hartshorne provides a counter-example to this conjecture in
[{\bf 5}]. He defined a module $N$ to be $\fa$-cofinite if
$\Supp_RN\subseteq \V(\fa)$ and $\Ext^i_R(R/\fa,N)$ is finitely
generated for all $i\geq 0$ and he asked the following question.

\begin{question} Let $\fa$ be an ideal of $R$ and $N$ a finitely
generated $R$-module. When are $H^i_{\fa}(N)$ $\fa$-cofinite for
all $i\geq 0$?
\end{question}

Hartshorne proved that if $\fa$ is an ideal of the complete
regular local ring $R$ and $N$ a finitely generated $R$-module,
then
$H^i_{\fa}(N)$ is $\fa$-cofinite in two cases: \\
(i) (see [{\bf 5}, Corollary 6.3]) $\fa$ is a principal ideal, and \\
(ii) (see [{\bf 5}, Corollary 7.7]) $\fa$ is a prime ideal with $\dim R/\fa=1$.\\

This subject was studied by several authors afterward ( see e.g.
[{\bf 8}],[{\bf 4}] and [{\bf 15}]). The best result concerning
cofiniteness of local cohomology is:

\begin{theorem} ([{\bf 8}],[{\bf 4}] and [{\bf 15}]) Let $\fa$ be an
ideal of $R$ and $N$ a finitely generated $R$-module. If either
$\fa$ is principal or $R$ is local and $\dim R/\fa=1$, then
$H^i_{\fa}(N)$ is $\fa$-cofinite for all $i\geq 0$.
\end{theorem}

S. Yassemi [{\bf 13}, Question 2.7] asked wether 1.2 holds for
generalized local cohomology. The main aim of this paper is to
extend 1.2 to generalized local cohomology. More precisely, we
prove the following.

\begin{theorem} Let $\fa$ denote an ideal of the  ring
$R$. Let $M$ and $N$ be two finitely generated $R$-modules with
$\pd M < \infty$. If either
\begin{itemize}
\item[(i)] $\fa$ is principle, or
\item[(ii)] $R$ is
complete local and $\fa$ is a prime ideal with $\dim R/\fa=1$,
\end{itemize}
then $H^i_{\fa}(M,N)$ is $\fa$-cofinite for all $i\geq 0$.
\end{theorem}

All rings considered in this paper are assumed to be commutative
Noetherian with identity. In our terminology we follow that of the
text book [{\bf 3}].

\section {Cofiniteness results}

Let $\fa$ denote an ideal of a ring $R$. The generalized local
cohomology defined by
$$H^i_{\fa}(M,N)=\underset{n}{\varinjlim}\Ext_R^i(M/\fa^nM,N),$$
for all $R$-modules $M$ and $N$. Note that this is in fact a
generalization of the usual local cohomology, because if $M=R$,
then $H^i_{\fa}(R,N)=H^i_{\fa}(N)$.

\begin{definition}  Let $M$ be an $R$-module. The generalized ideal
transform functor with respect to an ideal $\fa$ of $R$ is defined
by
$$D_{\fa}(M,\cdot)=\underset{n}{\varinjlim}\Hom_R(\fa^nM,\cdot).$$
\end{definition}

Let $R^iD_{\fa}(M,\cdot)$ denote the $i$-th right derived functor
of $D_{\fa}(M,\cdot)$. One can check easily that there is a
natural isomorphism $R^iD_a(M,\cdot)\cong
\underset{n}{\varinjlim}\Ext^i_R(\fa^nM,\cdot)$. Thus, by
considering the Ext long exact sequences induced by the short
exact sequences
$$0\lo \fa^nM\lo M\lo M/\fa^nM\lo 0, \quad (n\in\BN),$$ we  can
deduce the following lemma.

\begin{lemma}  Let $M$ be an $R$-module. For any $R$-module $N$, there is an
exact sequence
\begin{multline*}
0\lo H^0_{\fa}(M,N) \lo \Hom_R(M,N)\lo D_{\fa}(M,N)\lo
H^1_{\fa}(M,N)\lo\dots\\ \lo\dots\lo H^i_{\fa}(M,N)\lo
\Ext_R^i(M,N)\lo R^iD_{\fa}(M,N)\lo H^{i+1}_{\fa}(M,N)\lo \dots.
\end{multline*}
Moreover, if $M$ has finite projective dimension, then there is a
natural isomorphism $H^{i+1}_{\fa}(M,N) \cong R^iD_{\fa}(M,N)$ for
all $i\geq \pd M+1$.
\end{lemma}

Let $\fa$ be an ideal of $R$ and $M$ an $R$-module. We say that
$M$ is $\fa$-cofinite, if $\Supp_RM\subseteq V(\fa)$ and
$\Ext_R^i(R/\fa,M)$ is finitely generated for all $i\geq 0$.

\begin{lemma} Suppose $M,N$ are two $R$-modules and $\fa$ an ideal
of $R$. If $M$ is finitely generated, then
$\Supp_RH^i_{\fa}(M,N)\subseteq V(\fa)$, for all $i\geq 0$.
\end{lemma}

{\bf Proof.} Let $\fp$ be a prime ideal of $R$. It follows from
[{\bf 9}, Theorem 9.50], that $$\Ext_R^i(M/\fa^nM,N)_{\fp}\cong
\Ext_{R_{\fp}}^i(M_{\fp}/(\fa^nR_{\fp})M_{\fp},N_{\fp}),$$ for all
$i\geq 0$. On the other hand, it is well known that the formation
tensor product preserves direct limits (see e.g. [{\bf 9},
Corollary 2.20]). Thus
$$H^i_{\fa}(M,N)_{\fp}\cong R_{\fp}\otimes_R\underset{n}{\varinjlim}
\Ext_R^i(M/\fa^nM,N)\cong
\underset{n}{\varinjlim}\Ext_{R_{\fp}}^i(M_{\fp}/(\fa^nR_{\fp})M_{\fp},N_{\fp})
\cong H^i_{\fa R_{\fp}}(M_{\fp},N_{\fp}).$$ This shows that
$\Supp_RH^i_{\fa}(M,N)\subseteq V(\fa)$, as required.
$\blacksquare$

\begin{lemma}
\begin{itemize}
\item[(i)] If $M$ is a finitely generated $R$-module such that
$\Supp_RM\subseteq V(\fa)$, then $M$ is $\fa$-cofinite.
\item[(ii)] Let $0\lo L\lo M\lo N\lo 0$ be an exact sequences of $R$-modules.
Whenever two of $L,M$ or $N$ are $\fa$-cofinite, then the third
one is also $\fa$-cofinite.
\end{itemize}
\end{lemma}

{\bf Proof.} (i) Since $M$ is finitely generated, it follows that
$\Ext_R^i(R/\fa,M)$, is finitely generated for all $i\geq 0$.
Hence $M$ is $\fa$-cofinite, by definition.

(ii) This is well known and can be deduced easily, by considering
the long exact sequence $$ \dots\lo \Ext_R^i(R/\fa,L)\lo
\Ext_R^i(R/\fa,M)\lo \Ext_R^i(R/\fa,N)\lo \Ext_R^{i+1}(R/\fa,L)\lo
\dots .$$ $\blacksquare$

\begin{lemma} Let $\fa=Ra$ be a principal ideal of $R$ and $M$ and
$N$ two finitely generated $R$-modules. Let $\Hom_R(M,N)_a$ denote
the localization of $\Hom_R(M,N)$ with respect to the
multiplicative closed subset $\{a^i: i\geq 0\}$ of $R$. Then
\begin{itemize}
\item[(i)] there is a natural isomorphism $D_{\fa}(M,N)\cong
\Hom_R(M,N)_a$, and
\item[(ii)] $H^1_{\fa}(M,N)$ is
$\fa$-cofinite.
\end{itemize}
\end{lemma}

{\bf Proof.} (i) If $a$ is nilpotent, then it is clear that both
$D_{\fa}(M,N)$ and $\Hom_R(M,N)_a$ will vanish. Hence, we may and
do assume that $a$ is not nilpotent. For all $i,j\in \mathbb{N}$
with $j\geq i$, let $\pi_{ij}:\Hom_R(a^iM,N)\lo \Hom_R(a^jM,N)$ be
the map defined by $\pi_{ij}(f)=f|_{a^jM}$, for all $f\in
\Hom_R(a^iM,N)$. Also, denote the natural map $\Hom_R(a^iM,N)\lo
D_{\fa}(M,N)$, by $\pi_i$. Recall that, we defined $D_{\fa}(M,N)$
as the direct limit of the direct system
$(\Hom_R(a^iM,N),\pi_{ij})_{i,j\in \mathbb{N}}$.

Now define $\psi_i:\Hom_R(a^iM,N)\lo (\Hom_R(M,N))_a$, by
$\psi_i(f)=f\lambda_i/a^i$, where $\lambda_i:M\lo a^iM$ is defined
by $\lambda_i(m)=a^im$, for all $m\in M$. Clearly
$\{\psi_i\}_{i\in \mathbb{N}}$ is a morphism between direct
systems. Assume $\psi:D_{\fa}(M,N)\lo (\Hom_R(M,N))_a$ is the
homomorphism induced by $\{\psi_i \}_{i\in \mathbb{N}}$. Thus for
each $g\in D_{\fa}(M,N)$, we have $\psi (g)=\psi_i(f)$, where
$i\in \mathbb{N}$ and $f\in \Hom_R(a^iM,N)$ are such that
$\pi_i(f)=g$. We show that $\psi$ is an isomorphism. First, we
show that $\psi$ is injective. Suppose $\psi (g)=0$, for some
$g\in D_{\fa}(M,N)$. There are $i\in \mathbb{N}$ and $f\in
\Hom_R(a^iM,N)$, such that $g=\pi_i(f)$. Hence $$\psi
(g)=\psi_i(f)=f\lambda_i/a^i=0.$$ Hence there is $t\in \mathbb{N}$
such that $a^t(f\lambda_i)=0$. Set $j=i+t$. Then it follows that
$\pi_{ij}(f)=0$ and so
$$g=\pi_i(f)=\pi_j(\pi_{ij}(f))=0.$$

Next, we show that $\psi$ is surjective. Let $x_1,x_2,\dots ,x_t$
be a set of generators of $M$. Let $l\in (\Hom_R(M,N))_a$. Then
there are $h\in \Hom_R(M,N)$ and $c\in \mathbb{N}$, such that
$l=h/a^c$. Since $N$ is a Noetherian $R$-module, there exists an
integer $e\geq c$, such that $(0:_Na^e)=(0:_Na^{e+j})$, for all
$j\geq 0$. Define $f\in \Hom_R(a^{2e}M,N)$, by
$f(a^{2e}x)=a^{2e-c}h(x)$, for all $x\in M$. If
$a^{2e}x=a^{2e}x'$, for some $x$ and $x'$ in $M$, then $h(x-x')\in
(0:_Na^{2e})$. Hence $a^{2e-c}h(x)=a^{2e-c}h(x')$. Therefore $f$
is well-defined. Set $g=\pi_{2e}(f)$. Then
$$\psi (g)=\psi_{2e}(f)=f\lambda_{2e}/a^{2e}=h/a^c=l.$$
Thus $\psi $ is surjective.

(ii) Let $\psi:D_{\fa}(M,N)\lo (\Hom_R(M,N))_a$ be as above.
 By part (i), [{\bf 3}, Theorem 2.2.4(i)] and
2.2 we have the following commutative diagram with exact rows.
\begin{equation*}
\setcounter{MaxMatrixCols}{11}
\begin{smallmatrix}
0 &\lo &H^0_{\fa}(M,N) &\lo &\Hom_R(M,N)&\overset{f}{\lo}&D_a(M,N)
&\lo & H^1_{\fa}(M,N)& \stackrel{g}{\lo} &\Ext_R^1(M,N)
\\ & & & &\mapdown{id} & &\mapdown{\psi} & & & & & &
\\ 0 &\lo &\Gamma_{\fa}(\Hom_R(M,N)) &\lo & \Hom_R(M,N)
&\stackrel{h}{\lo}
&\Hom_R(M,N)_a & \lo& H^1_{\fa}(\Hom_R(M,N))& \lo 0. &
\end{smallmatrix}
\end{equation*}

Let $K$ be the kernel of the map $g$. We have $K\cong \coker f$
and $H^1_{\fa}(\Hom_R(M,N))\cong \coker h$.  The map $\psi$
induces an isomorphism ${\psi}^*:\coker f\lo \coker h$, which is
defined by ${\psi}^*(x+im f)=\psi(x)+im h$, for all $x+im f\in
\coker f$. Hence, it follows that $K\cong H^1_{\fa}(\Hom_R(M,N))$.
Therefore $K$ is $\fa$-cofinite, by 1.2. Now consider the exact
sequence
$$0\lo K\lo H^1_{\fa}(M,N)\lo \im g\lo 0.$$ Since $\Ext_R^1(M,N)$
is finitely generated, it follows by, 2.3 and 2.4(i), that $\im g$
is $\fa$-cofinite. Thus $H^1_{\fa}(M,N)$ is $\fa$-cofinite, by
2.4(ii). $\blacksquare$

\begin{lemma}  Let $\fa$ denote an ideal of the  ring
$R$ and $N$ an $\fa$-cofinite $R$-module. Suppose that for any
finitely generated $R$-module $M$ with $\pd M<\infty$,
$\Hom_R(M,N)$ (resp. $M\otimes_RN)$ is $\fa$-cofinite.
 Then $\Ext_R^i(M,N)$ (resp. $\Tor_i^R(M,N)$) is $\fa$-cofinite for
all finitely generated $R$-modules $M$ with $\pd M<\infty$ and all
$i\geq 0$.
\end{lemma}

{\bf Proof.} We prove only the $\fa$-cofiniteness of
$\Ext_R^i(M,N)$, $i\geq 0$ and the proof of the other part is
similar. The proof proceeds by induction on $t=\pd M$. For $t=0$,
the claim holds by assumption. Now, suppose $t>0$. There is a
short exact sequence $$0\lo K\lo R^n\lo M\lo 0.$$ From this
sequence, we deduce the exact sequence $$0\lo \Hom_R(M,N)\lo
\Hom_R(R^n,N)\lo\Hom_R(K,N)\lo \Ext_R^1(M,N)\lo 0,$$ and the
isomorphisms $\Ext_R^{i+1}(M,N)\cong\Ext_R^i(K,N)$ for all $i\geq
1$. Thus from induction hypothesis, we deduce that
$\Ext_R^{i+1}(M,N)$ is $\fa$-cofinite for all $i\geq 1$. Note that
$\pd K<t$. Also, by using the above exact sequence, one can check
easily that $\Ext_R^1(M,N)$ is $\fa$-cofinite. Therefore, the
claim follows by induction. $\blacksquare$

\begin{lemma} Let $\fa$ denote an ideal of the ring
$R$. Let $M$ and $N$ be two finitely generated $R$-modules with
$\pd M<\infty$. If either
\begin{itemize}
\item[(i)] $\fa$ is principal, or
\item[(ii)]$R$ is
complete local and $\fa$ is a prime ideal with $\dim R/\fa=1$,
\end{itemize}
then $\Ext_R^p(M,H^q_{\fa}(N))$ is $\fa$-cofinite for all $p,q\geq
0$
\end{lemma}

{\bf Proof.} First, we consider the case that $\fa$ is principal.
By [{\bf 9}, Theorem 11.38], there is a Grothendieck spectral
sequence
$$E_2^{p,q}:=\Ext_R^p(M,H^q_{\fa}(N))\underset{p}{\Longrightarrow}
H_{\fa}^{p+q}(M,N).$$  We have $E_2^{p,q}=0$ for $q\neq 0,1$,
because, by [{\bf 3}, Theorem 3.3.1], $H_{\fa}^q(N)=0$, for all
$q>1$. Since $E_2^{p,0}$ is finitely generated and
$\Supp_RE_2^{p,0}\subseteq V(\fa)$, it follows by 2.4(i), that
$E_2^{p,0}$ is $\fa$-cofinite. Therefore, it is enough to show
that $E_2^{p,1}$ is $\fa$-cofinite for all $p\geq 0$. By [{\bf 9},
Corollary 11.44], we have an exact sequence $$0\lo E_2^{1,0} \lo
H_{\fa}^1(M,N)\overset{g}{\lo} E_2^{0,1}\lo E_2^{2,0}
\overset{f}{\lo} H_{\fa}^2(M,N).$$ By 2.5(ii), $H_{\fa}^1(M,N)$ is
$\fa$-cofinite. Thus, it turns out that $\im g$ is $\fa$-cofinite,
by 2.4(ii). From the exact sequence
$$0\lo\im g\lo E_2^{0,1}\lo \ker f\lo 0,$$ we deduce that
$E_2^{0,1}$ is $\fa$-cofinite. Note that $\ker f$ is a finitely
generated $R$-module. Therefore 2.6 implies that $E_2^{p,1}$ is
$\fa$-cofinite for all $p\geq 0,$ because $H^1_{\fa}(N)$ is
$\fa$-cofinite by 1.2.

Now suppose that $R$ is a complete local ring and $\fa$ a prime
ideal of $R$ with $\dim R/\fa=1$. In view of 2.6, it suffices to
show that $\Hom_R(M,H^q_{\fa}(N))$ is $\fa$-cofinite for all
finitely generated $R$-modules $M$ with $\pd M<\infty$. We prove
this claim by induction on $\pd M=t$. The case $t=0$, is clear by
1.2. Now assume that $t>0$ and consider the exact sequence
$$0\lo K\lo R^n\lo M\lo 0.$$ It follows that $\pd K \leq t-1$. This
short exact sequence yields the exact sequence
$$0\lo\Hom_R(M,H^q_{\fa}(N))
\lo\Hom_R(R^n,H_{\fa}^q(N))\overset{f}\lo \Hom_R(K,H_{\fa}^q(N))
.$$ Since, by [{\bf 4}, Theorem 2] the subcategory of
$\fa$-cofinite $R$-modules is abelian, it follows that $\ker
f\cong \Hom_R(M,H^q_{\fa}(N))$ is $\fa$-cofinite. Therefore the
claim follows by induction. $\blacksquare$

\begin{theorem} Let $\fa$ denote a principal ideal of the ring $R$.
Let $M$ and $N$ be two finitely generated $R$-modules with $\pd
M<\infty$. Then $H^p_{\fa}(M,N)$ is $\fa$-cofinite for all $p\geq
0$.
\end{theorem}

{\bf Proof.} By [{\bf 9}, Theorem 11.38], there is a Grothendieck
spectral sequence
$$E_2^{p,q}:=\Ext_R^p(M,H^q_{\fa}(N))\underset{p}{\Longrightarrow}
H_{\fa}^{p+q}(M,N)$$ This implies the following exact sequence  in
view of [{\bf 11}, Ex. 5.2.2]. Note  that $E_2^{p,q}=0$ for $q\neq
0,1.$ $${\lo} E_2^{p,0} \overset{f}{\lo} H_{\fa}^p(M,N)
\overset{d}{\lo} E_2^{p-1,1}{\lo} E_2^{p+1,0} \overset{g} \lo
H_{\fa}^{p+1}(M,N)\lo\dots .$$ Now, $\im f$ is a quotient of
$E_2^{p,0}$ and so  is finitely generated. Hence $\im f$ is
$\fa$-cofinite, by 2.3 and 2.4(i). Also, $\ker g$ is
$\fa$-cofinite by the same reason. By considering the short exact
sequence $$0\lo \im d\lo E_2^{p-1,1} \lo \ker g\lo 0,$$ we deduce
that $\im d$ is $\fa$-cofinite. Note that $E_2^{p-1,1}$ is
$\fa$-cofinite by 2.7(i). Now from the short exact sequence
$$0\lo\im f\lo H_{\fa}^p(M,N) \lo\im d\lo 0,$$ we deduce that
$H^p_{\fa}(M,N)$ is $a$-cofinite for all $p\geq 0$. $\blacksquare$

\begin{theorem} Let $\fp$ denote a prime ideal of the complete local ring $(R,\fm)$
with $\dim R/\fp=1$, and $M,N$ two finitely generated $R$-modules
with $\pd M<\infty$. Then $H^i_{\fp}(M,N)$ is $\fp$-cofinite for
all $i\geq 0$.
\end{theorem}

{\bf Proof.} There is a spectral sequence
$$E_2^{p,q}:=\Ext_R^p(M, H_{\fp}^q(N)) \underset{p}{\Longrightarrow}
H^{p+q}_{\fp}(M,N)=E^n.$$ It follows from 2.7(ii) that $E_2^{p,q}$
is $\fp$-cofinite for all $p,q$. By considering the sequence
$$\dots\lo E_2^{p-2,q+1} \overset{d_2^{p-2,q+1}}{\lo} E_2^{p,q}
\overset{d_2^{p,q}} \lo E_2^{p+2,q-1}\lo\dots ,$$ we deduce that
$\im d_2^{p-2,q+1}$ and $\ker  d_2^{p,q}$ are $\fp$-cofinite, by
[{\bf 4}, Theorem 2]. Hence $E_3^{p,q}=\ker d_2^{p,q}/\im
d_2^{p-2,q+1}$ is $\fp$-cofinite. By irritating this arguments we
get that $E_r^{p,q}=\ker d_{r-1}^{p,q}/\im d_{r-1}^{p-r+1,q+r-2}$
is $\fp$-cofinite for all $r>0$ and so $E_\infty^{p,q}$ is
$\fp$-cofinite for all $p,q\geq 0$. There is a filtration
$$E^n=E_0^n\supseteq \dots \supseteq E_p^n \supseteq \dots
\supseteq E_n^n\supseteq E_{n+1}^n =0,$$ such that
$E_p^n/E_{p+1}^n\cong E_\infty^{p,n-p}$.  Thus $E_n^n$ is
$\fp$-cofinite. Now, by applying 2.4(ii) repeatedly on the short
exact sequences $$0\lo E_{p+1}^n\lo E_p^n\lo E_{\infty}^{p,n-p}\lo
0,   p=0,1\dots ,n-1,$$ we deduce that $E^n$ is $\fp$-cofinite, as
required. $\blacksquare$

Many  results concerning local cohomology in positive prime
characteristic can be extend to generalized local cohomology. In
particular the main results of [{\bf 7}] are also hold for
generalized local cohomology.

\begin{theorem} Let $(R,\fm)$ be a regular
local ring of characteristic $p>0$, and $\fa$  an ideal of $R$.
Let $\fp$ be a prime ideal of $R$ and $M$ a finitely generated
$R$-module. Then
\begin{itemize}
\item[(i)]
$\mu^i(\fp,H_{\fa}^j(M,R))\leq\mu^i(\fp,\Ext_R^j(M/\fa
 M,R))$ for all $j\geq 0$. In particular $\mu^i(\fp,H^j_{\fa}(M,R))$ is
finite for all $j\geq 0$ and all $i\geq 0$.
\item[(ii)]
$\Ass_R(H^j_{\fa}(M,R))\subseteq\Ass_R(\Ext_R^j(M/\fa M,R))$ and
so $\Ass_R(H_{\fa}^j(M,R))$ is finite for all $j\geq 0$.
\end{itemize}
\end{theorem}

{\bf Proof.} The proof is a straightforward adoption of the proof
of [{\bf 7}, Theorem 2.1 and Corollary 2.3]. $\blacksquare$

{\bf Acknowledgments.} We thank the referee whose comments make
the paper more readable. Also, the first named author would like
to thank the Institute for Studies in Theoretical Physics and
Mathematics (Tehran) for financial support.


\end{document}